\documentclass[11pt]{amsart}

\usepackage{amssymb}
\usepackage{amsmath}
\usepackage{amsfonts}

\input{epsf}
\usepackage{graphicx}
\usepackage{graphics}

\usepackage{psfrag}
\usepackage{epsfig}

\makeatletter 
 
 \@addtoreset{equation}{section}
\makeatother

\newtheorem{theorem}{Theorem}[section]
\newtheorem{lemma}[theorem]{Lemma}

\newtheorem{cor}[theorem]{Corollary}

\theoremstyle{definition}

\newtheorem{remark}[theorem]{Remark}
\newtheorem{example}[theorem]{Example}

\newcommand{\HFK}{\widehat{HFK}}
\newcommand{\HF}{\widehat{HF}}
\newcommand{\CF}{\widehat{CF}}
\newcommand{\on}{\operatorname}

\newcommand{\Spinc}{\on{Spin}^c}

\newcommand{\um}{{\bf v}_-}
\newcommand{\up}{{\bf v}_+}

\newcommand{\ZZ}{\mathbb{Z}}

\newcommand{\QQ}{\mathbb{Q}}
\newcommand\goth[1]{\mathfrak{#1}}
\newcommand{\s}{\goth{s}}

\newcommand{\sign}{\on{sign}}

\newcommand{\rk}{\on{rk}}

\begin{document}

\title{Khovanov homology and tight contact structures}

\author{Olga Plamenevskaya}
\address{Department of Mathematics, SUNY Stony Brook, Stony Brook, NY 11794}
\email{olga@math.sunysb.edu}

\begin{abstract} Using the relation between Khovanov homology and the Heegaard Floer homology of 
branched double covers, we show how Khovanov homology can be used to establish tightness of branched 
double covers of certain transverse knots. We give examples of several infinite families of knots 
whose branched covers are tight for Khovanov-homological reasons, and show that some of these branched 
covers are not Stein fillable. 

\end{abstract}

\maketitle

\section{Introduction}

The goal of this paper is to demonstrate how  Khovanov homology can be used to establish  tightness 
of certain contact structures. The contact manifolds we study are all branched double covers 
of transverse knots in standard contact $S^3$; our main tool is the relation between the Khovanov-homological 
invariant $\psi$ \cite{Pla1} and the Heegaard Floer contact invariant $c(\xi)$ \cite{ContOS}
of the branched double cover of $K$. This relation was conjectured by the author in 
\cite{Pla2} and proved by Lawrence Roberts \cite{Ro}; 
the following theorem is a corollary of the results in \cite{Ro}. 

\begin{theorem} Suppose that the spectral sequence \cite{OS2} from the reduced Khovanov homology 
to the Heegaard Floer homology $\HF(-\Sigma(K))$ of the branched double cover of $K$ collapses 
at the $E^2$ stage, thus providing an isomorphism between the Khovanov and the Heegaard Floer homology.
(In particular, this is true when $K$ belongs to a quasi-alternating knot type \cite{OS2}.) 
Then $c(\xi)\neq 0$ whenever $\psi\neq 0$.   \label{roberts} 
\end{theorem}

We work with $\ZZ/2$ coefficients throughout; $\psi$ denotes the version of the invariant  from \cite{Pla1}
that lives in the reduced Khovanov homology with $\ZZ/2$ coefficients. 

We will use Theorem \ref{roberts} together with the following non-vanishing 
criterion for $\psi$: 

\begin{theorem} \label{nonzero} 
If $K$ is a transverse knot such that $sl(K)=s-1$,  then $\psi(K) \neq 0$. The converse is 
also true if $K$ belongs to a quasi-alternating (or any $Kh_{\ZZ/2}$-thin) smooth knot type.
Here $s=s(K)$ stands for Rasmussen's invariant \cite{Ra}. 
\end{theorem}

\begin{cor} \label{tightness} If $K$ is a transverse representative of a quasi-alternating knot, the induced contact 
structure $\xi_K$ on the branched  double cover of $K$ has $c(\xi) \neq 0$ if  $sl(K)=\sigma-1$,
where $\sigma$ is the signature of the knot (with the sign conventions such that $\sigma(\text{right trefoil})=2$). 
In particular, $\xi_K$ is tight whenever $sl(K)=\sigma-1$.  
\end{cor}

Note that for any transverse knot  $sl(K) \leq s-1$ \cite{Pla1, Sh}, 
so the condition of Theorem~\ref{nonzero}
is equivalent to sharpness of this upper bound for the self-linking number.  

The values of the maximal self-linking number for knots with 10 crossings or less are known   
\cite{Ng}; combining those with the above results, we can in certain cases establish the existence 
of a tight contact structure on the branched double cover of a given smooth knot.

To demonstrate the efficiency of Khovanov homology in proofs of tightness, we need to find examples of transverse knots satisfying the hypothesis of Corollary \ref{tightness} (or those of Theorems \ref{roberts} and \ref{nonzero}). 
Many of these are provided by quasipositive braids, but the corresponding contact structures are obviously 
Stein fillable and therefore tight. However, non-quasipositive braids with $sl(K)=s-1$ do exist; 
in section \ref{examples}, we give examples of several infinite families of such transverse knots. 
We also show that some of the corresponding contact structures are not Stein fillable (and thus the result 
is non-trivial). Some of our examples are transverse 3-braids; accordingly, their branched covers have open 
book decomposition of genus one whose tightness can be established 
by methods \cite{Ba, HKM}. Unlike \cite{Ba, HKM}, our proofs work 
for transverse braids of arbitary index,  require no explicit calculations of the Heegaard Floer contact invariants, and are completely combinatorial (once Theorem \ref{roberts}
is in our hands). In fact, John Baldwin has recently 
written a computer program \cite{Ba1} for determining whether $\psi$ is non-zero for a given transverse 
braid. If the underlying link is quasi-alternating (or otherwise satisfies the collapsing condition of Theorem \ref{roberts}), we can establish tightness of the corresponding contact structure by a computer calculation.    

\smallskip

\noindent {\bf Acknowledgements.} I would like to thank John Baldwin, Lenny Ng and Andras Stipsicz for 
very helpful conversations, and Lawrence Roberts for helpful correspondence.    

\section{Non-vanishing of the transverse and contact invariants}

Theorem \ref{roberts} is essentially a corollary of Roberts's work \cite{Ro}. However, it is not contained 
in Roberts's paper, and as the  constructions of \cite{Ro} require some care, we find it useful to 
review the results of 
\cite{Ro} before giving  a proof of Theorem \ref{roberts}. Our notation below  is a bit different from that of \cite{Ro};
we will also ignore a few  minor details (e.g. for technical 
reasons one needs to add two extra strands to a given transverse braid, etc).     
 
We first recall the construction of spectral sequence
that relates the reduced Khovanov homology and the Heegaard Floer homology of the branched double cover \cite{OS2}.
The main result of \cite{OS2} gives a link surgeries filtered chain complex $(C(K), D)$ whose associated spectral sequence  converges  to
$\HF(-\Sigma(K))$ and has the $E^1$-term
$$
E^1 = \bigoplus_{i \in \{0,1\}^n} \HF(-\Sigma(K_i)), 
$$
where the sum is taken over all complete resolutions $K_i$ of the knot (or link) $K$. (As usual, the components 
of $i = (i_1, \dots, i_n)$ stand for $0$- and $1$- resolutions of crossings of $K$.) The filtration of 
$C(K)$ is given by the flattened cube filtration $I = \sum_{j=1}^n {i_j}$; the differential 
does not decrease this $I$-filtration. The complex $(C(K), D)$ is constructed by counting holomorphic polygons
in  (a symmetric product of) a Heegaard diagram compatible with all link surgeries; as a vector space, 
$$
C(K) = \bigoplus_{i \in \{0,1\}^n} \CF(-\Sigma(K_i)), 
$$
and the differential $D^0$ on the associated graded object is given by the sum of the usual Heegaard Floer 
boundary maps $d: \CF(-\Sigma(K_i)) \to \CF(-\Sigma(K_i))$. 
One  shows that for each complete resolution of the knot $K$, $\HF(-\Sigma(K_i))$ is precisely 
the component $CKh(K_i)$ of the reduced Khovanov complex. (In fact, each manifold $-\Sigma(K_i)$
is simply the connected sum of a few copies of $S^1\times S^2$.) The differential $D_1$ on   
$\bigoplus_{i \in \{0,1\}^n} \HF(-\Sigma(K_i))$ given by sums of maps 
$$
\HF(-\Sigma(K_i)) \to  \HF(-\Sigma(K_i'))
$$ with $I(i') = I(i)+1$; each such  map on $\HF$ is induced by 
the surgery cobordism between  $\HF(-\Sigma(K_i))$ and $\HF(-\Sigma(K_i'))$ that corresponds to the change of 
resolution relating $K_i$ and $K_i'$. Under the correspondence $\HF(-\Sigma(K_i))= CKh(K_i)$, these maps 
$\HF(-\Sigma(K_i)) \to \HF(-\Sigma(K_i'))$ are shown to be equal to maps $CKh(K_i) \to CKh(K_i')$ that 
form the differential on the reduced Khovanov complex. This implies that the $E^1$-term of the complex $C(K)$
is the reduced Khovanov complex, and the $E^2$-term is the reduced Khovanov homology. 
 
The main idea of \cite{Ro} is to endow the above complex $C(K)$ with an additional filtration induced by the 
binding of an open book decomposition of $\HF(-\Sigma(K))$. More precisely, represent the knot $K$ as 
a (transverse) braid in $S^3$,  and let $B$ be the braid axis. We can always 
assume (stabilizing if necessary) that the braid index of $K$ is odd. Then the manifold $\HF(-\Sigma(K))$ has a natural open 
book decomposition with binding $B$ and pages given by branched double covers of a disk; this open book is compatible with the contact structure on $\Sigma(K)$ induced by $K$. One can incorporate the knot $B$ into constuctions 
of \cite{OS2}, i.e. consider link surgeries  Heegaard diagrams compatible with $B$.  Then the knot $B$ induces the additional ``Alexander'' filtration 
on $C(K)$. 
In particular,  $B$ filters $E_I^1= \bigoplus_{i \in \{0,1\}^n} \HF(-\Sigma(K_i))$, and the cobordism 
maps   $\HF(-\Sigma(K_i)) \to \HF(-\Sigma(K_i'))$ considered above respect the filtration. 
(It is important to note here that  $\HF(-\Sigma(K_i)) = \HFK(-\Sigma(K_i), B)$, 
i.e. the knot homology is trivial and only gives a filtration of $\HF(-\Sigma(K_i))$ in this case.)
We keep the notation $C(K)$ for the bi-filtered complex, and use the subscripts $A$ and $I$ to distinguish 
between spectral sequences associated to different filtrations.

Roberts shows that as a bi-filtered complex, the term $E_I^1$ is isomorphic (modulo some  grading adjustments)
to the bi-filtered ``skein Khovanov complex''  \cite{APS},  which is the usual reduced Khovanov complex  
$\bigoplus_{i \in \{0,1\}^n} CKh(K_i)$ endowed with an extra filtration.
The extra filtration is induced by $B$ and comes from dividing the components of each resolution 
$K_i$ of $K$ into those circles that link with $B$ and those that do not. 
It turns out that for each $i$, this filtration coincides  with 
the Alexander filtration on knot Floer homology  $\HF(-\Sigma(K_i))= \HFK(-\Sigma(K_i), B)$, thus 
the skein Khovanov complex is indeed the $E_I^1$-term for $C(K)$. 
 
As before, we can consider the spectral sequence induced on  $C(K)$ by the $I$-filtration; we still 
have $E_I^2= Kh(K)$, and $E_I^\infty = \HF(-\Sigma(K))$, but now each page of the spectral sequence inherits a filtration from the $A$-filtration on the original complex. (For a filtered chain complex, a filtration 
on its homology group is defined, as usual, by taking the minimum of filtration levels of cycles representing 
a given homology class.) We point out a caveat here: the $A$-filtration does not behave well with respect 
to the spectral sequence; in particular, when $E_I^n = E_I^\infty = \HF(-\Sigma(K))$, the $A$-filtration 
on $E_I^n$ computed from the spectral sequence may differ from the $A$-filtration on $\HF(-\Sigma(K))$ computed
by taking the homology of the entire complex (see Remark \ref{wrong} below).  

On the other hand, we can ignore the $I$-filtration and consider the spectral 
sequence from  
$C(K)$ to 
$\HF(-\Sigma(K))$ 
induced by the $A$-filtration. Then the $E_A^1$-term is given by  $\HFK(-\Sigma(K), B)$, since  
 the link surgeries construction of \cite{OS2} work just as well when we pass to 
the associated graded object for the knot filtration of $B$ \cite[Proposition 7.1]{Ro}. 
Moreover, the subsequent pages of this spectral sequence are quasi-isomorphic to the pages of the knot Floer homology 
spectral sequence induced by the knot filtration of $B$ on the Heegaard Floer complex  $\CF(-\Sigma(K))$ 
\cite[Lemma 7]{Ro}.

Now, recall that the contact invariant $c(\xi)$ is the image in $\HF(-\Sigma(K))$ of the unique lowest $A$-filtration 
element $c \in \HFK(-\Sigma(K), B)$ (which lies in the $A=-g(B)$ filtration level).  Accordingly, the $A=-g(B)$ filtration level in $\HF(-\Sigma(K))$
is empty or one-dimensional depending on whether $c(\xi)$ vanishes or not. In the $A$-filtered skein $CKh(K)$, we can 
also pinpoint the unique lowest $A$-degree generator $\psi \in \bigoplus_{i \in \{0,1\}^n} CKh(K_i)$.  
Indeed, the construction of the skein filtration on  $CKh(K)$ implies that the lowest $A$-degree 
can be only reached when we take the oriented resolution of the braid $K$ 
(so that all resulting circles link with $B$), and pick the lowest 
quantum degree element $\psi = \um\otimes\dots\otimes \um$ in the corresponding component  $CKh(K_i)$.
Observe that this is precisely the cycle in the reduced Khovanov complex that gives the transverse invariant 
$\psi \in Kh(K)$ of \cite{Pla1}. (Strictly speaking, in the reduced case we take  $\psi = \um\otimes\dots\otimes \um \oplus \up$, with the $\up$ on the marked circle, but we may keep the notation without $\up$ by identifying 
the reduced complex with the subcomplex $CKh_{-}(K)$, see \cite{Pla1}  for details.) 
It then follows that in  the (A-filtered) $E^2_I=Kh(K)$-term of the spectral sequence induced on $C(K)$ by the $I$-filtration, 
the  $A=-g(B)$ filtration level is empty or one-dimensional depending 
on whether $\psi$ vanishes or not.






\begin{remark} \label{wrong} Even though both $\psi(K)\in E_I^2$ and $c(\xi_K) \in E_I^\infty$ lie 
in the same lowest $A$-filtration level and  are the images under the spectral sequence of the 
same canonical cycle in the skein Khovanov complex $E_I^1$, one should use caution when talking about 
the ``correspondence'' between $\psi$ and $c(\xi)$. Indeed, as was pointed out by John Baldwin, 
it is possible that $\psi(K)=0$ while $c(\xi_K)\neq 0$. This happens when $\psi$ is the boundary of a cycle $x$ in $E_I^1$, and $Dx$ has other terms of higher filtration level in the entire chain complex, so that $c(\xi)$ is not 
a boundary. A toy example of this phenomenon is given by a complex $(C, d)$ generated by three elements $x_{0,0}$, 
$y_{-2,1}$, $z_{-1, 2}$, where the indices indicate the $(A, I)$ bi-filtration, with $dx = y+z$. The lowest 
$A$-filtration element $y$ plays the role of $c(\xi)$. In the $I$-induced spectral sequence, the differential 
$d_0$ is trivial, $E_I^1$ is generated by classes $[x]$, $[y]$, $[z]$, and $d^1[x]=[y]$. Thus $[y]$ is 
a boundary in   $E_I^1$, while $y$ is not a boundary in the entire complex $(C, d)$. We also observe that 
even though the spectral sequence collapses at the $E^2$-term, $E_I^2$ and $H_*(C, d)$ are different as 
$A$-filtered vector spaces: the filtration level of the generator $[z]$ of $E_I^2$ is $A=-1$, while the 
filtration level of the generator of  $H_*(C, d)$ is $A=-2$.  

(Baldwin found an explicit family of transverse braids with $\psi(K)=0$ and $c(\xi_K)\neq 0$; we return 
to his examples in the next section.) However, we will show below that  when the spectral sequence for the $I$-filtration collapses at the $E^2$-term, non-vanishing of $\psi$ implies non-vanihing of $c(\xi)$.     
\end{remark}     

\begin{proof}[Proof of Theorem \ref{roberts}] The proof becomes easy if the chain complex $(C(K), D)$ satisfies the 
additional condition
\begin{equation}\label{cond}\tag{*}
\begin{split}  
&\text{the lowest $A$-filtration level  of  the complex $C(K)$ is $A=-g$;} \\  
&\text{there is a unique generator $c$ with $A(c) =-g$}. 
\end{split}
\end{equation}
If (\ref{cond}) holds,  the  element $c$ is  necessarily a $D$-cycle concentrated in a single $I$-grading level, 
and thus gives rise to a cycle in every term of the spectral 
sequence. Moreover, because $A(\psi)=-g$  for the  transverse element $\psi\in E^1_I$,    
$c$ is a representative of the class $\psi$, and $I(c)=0$.  

To ensure that (\ref{cond}) holds, we may need to pass to a new complex  $(C'(K), D')$, using a cancellation lemma 
(Lemma \ref{cancel} below). In fact, applying this lemma to  $(C(K), D)$ 
we obtain  the complex $C'(K)= \bigoplus_{i\in \{0,1\}^n} \HFK(-\Sigma(K_i), B)$. 
(Note that because  $\HF(-\Sigma(K_i))=\HFK(-\Sigma(K_i), B)$, the differential $(D')_I^0$ is trivial, so that 
the complex $C'(K)$ is the same as the $E_I^1$-term for $C(K)$.) 
It is now clear that $C'(K)$ has a unique lowest $A$-filtration element $c$; this element has $A=-g$ and lies in the component
$\HFK(-\Sigma(K_i), B)$  corresponding to the oriented resolution of the braid $K$. 

Condition (4)  of Lemma \ref{cancel} implies that $\psi \neq 0$ if and only if  
$E^2$-term of the $I$-induced spectral sequence on 
$(C'(K), D')$ is non-empty in the filtration level $A=-g$, and $c(\xi)\neq 0$ if and only if 
the  $E_I^\infty$-term of $(C'(K), D')$  is non-empty in the same filtration level, because the same is true for 
$(C(K), D)$. We will now assume 
that $(C(K), D)$ satisfies $(\ref{cond})$.

Now, suppose that $c(\xi)=0$. Then the $A=-g$ filtration level of $H_*(C(K), D))$ is empty, 
and thus the cycle $c$ is a boundary, $c=Dx$ for some $x\in C(K)$. 
Because $D$ is non-decreasing on the $I$-filtration, $I(x)\geq 0$.  We claim that $I(x)>0$: otherwise for the class $[x]$ in the associated graded object $E_I^0$ we have $D^0[x]= [Dx$ modulo terms of filtration level $I<0] = [c]$, which contradicts the fact 
that the class of $c$ (in $E_I^1$) is the non-zero cycle $\psi$. Thus $I(x)>0$, and $D^0[x]= 0$, which means 
that $[x]$ is a cycle that gives rise to an element of $E^1_I$. Next, consider $D^1[x]$. If $I(x)=1$, 
then  $D^1[x]= [Dx \text{ modulo terms of filtration level $I<0$}] = [c]=\psi$, which contradicts the hypothesis 
that the transverse element $\psi \in E_I^1$ is not a boundary in the Khovanov skein complex $(E_I^1, D^1)$. 
Thus we conclude that $I(x)\geq 2$, and  $D^1[x]=0$, so that the class of $[x]$ is a cycle in $E_I^1$ giving 
rise to an element of $E_I^2$. But since $Dx=c$, and $I(c) \leq I(x) -2$, it follows that the spectral sequence 
does not collapse at the $E_I^2$-term, a contradiction. 
     
\end{proof} 

To state the cancellation lemma used above, we consider a finitely generated bi-filtered complex $(C,d)$ with 
 an ascending filtration $I$ and a descending filtration $A$,  
\begin{eqnarray*}
&I: \quad  \ldots \subset C_{i+1} \subset C_{i} \subset \ldots \\
&A: \quad  \ldots \subset C_{a-1} \subset C_{a} \subset \ldots, 
\end{eqnarray*}
and let $C_{i,a}= (C_{i} \cap C_{a})/(C_{i+1} \cup C_{a-1})$ denote the bi-filtered quotients. Let 
$d_{00}:  C_{i,a} \to C_{i,a}$ stand for the differential induced by $d$. 

\begin{lemma} {\rm (\cite[Lemma 4.5]{Ra1}, \cite[Lemma 8]{Ro})} \label{cancel} Let $C$  be a 
finitely generated bi-filtered 
complex (over $\ZZ/2$).  There exists a unique (up to isomorphism) bi-filtered complex $C'$ such that 
\begin{enumerate}
\item $C$ and $C'$ are bi-filtered chain homotopy equivalent  
\item $(C')_{i,a}=H_*(C_{i,a})$ 
\item The differential $d'_{00}$ on the quotients of $C'$ is trivial. 
\item For each filtration, the spectral sequences
for $(C,d)$ and $(C', d')$ have the same terms, i.e. 
$$
E_I^r = (E'_I)^r \text{ and }E_A^r = (E'_A)^r \text{ for all } r\geq 0.  
$$
\end{enumerate}
\end{lemma}  
\qed

\begin{proof}[Proof of Theorem \ref{nonzero}] Suppose that $K$ is a transverse knot with $sl=s-1$. 
We are interested in the reduced Khovanov homology with $\ZZ/2$ coefficients, but it is convenient 
to consider the case of non-reduced homology with rational coeffients first. (Our notation will often be 
the same for different flavors of Khovanov homology; it will be clear for the context which case 
we are considering.) Recall that Lee  \cite{Lee}
introduces a differential $d'= d + \Phi$ on the Khovanov complex $CKh(K)$, where $d$ is the usual Khovanov's 
differential \cite{Kh}, and $\Phi$ is a map that raises the quantum grading. We recall that $d'$ corresponds 
to the multiplication and co-multiplication maps given by     
\begin{eqnarray*}
 m(\up\otimes \up) = m(\um\otimes \um)= \up \qquad    &\Delta(\up) = \up\otimes \um + \um\otimes\up \\ 
 m(\up\otimes \um) = m(\up\otimes \um)= \um \qquad    &\Delta(\um) = \um\otimes \um + \up\otimes\up 
\end{eqnarray*}
(we follow the notation from \cite{Ra} where $\um$ and $\up$ stand for the elementary generators 
of quantum degree $-1$ resp. $+1$.) 
This gives rise to a filtration in the Khovanov complex and yields a spectral sequence whose $E^2$-term is the 
usual $Kh(K)$, and the $E^\infty$ 
term is $Kh'=H_{*}(CKh, d')$. For a knot $K$, Lee's homology $Kh'(K)= \QQ \oplus \QQ$ is generated by two canonical cycles. Let $\s_0$ be a canonical generator corresponding 
to a choice of orientation of $K$; then $\psi$ is a $q$-homogeneous part of $\s_0$ with the lowest $q$-grading.
Indeed, the invariant $\psi$ is defined by representing the transverse knot by a braid, taking the oriented 
resolution, and taking the cycle
$$
\psi=\um \otimes \um \otimes \um \dots
$$
to be the lowest quantum degree term of the corresponding component of $CKh(K)$. 
The oriented resolution of a braid $K$ consists of nested circles, and 
$$
\s_0= (\um+\up) \otimes (\um-\up) \otimes (\um+\up) \dots
$$ 
is an element of $CKh$ obtained by labeling these circles by  $\um+\up$ and $\um-\up$, in alternating order. (The label on the outermost circle 
is determined by the orientation of the knot.)  Rasmussen \cite{Ra} defines a function $s$ on $Kh'$ whose value 
$s(x)$ on $x\in Kh'$ is the largest $n$ such that $x$ can be represented by a cycle all of 
whose terms have quantum grading at least $n$. The invariant $s$ is then defined so that 
$s-1=s_{min}=s([\s_0])$.  Recall  that $q(\psi)=sl(K)$ \cite{Pla1}, so our assumption means that 
$q(\psi)=s-1$. Suppose that $\psi(K)$ vanishes in $Kh(K)$, so $\psi = dy$ for some $y \in CKh(K)$. Since
$d$ preserves quantum gradings, we must have $q(y)=s-1$. Consider the element $\s_0-d'y$, where $d'$ is 
Lee's differential on $CKh$. This is a cycle in $Kh'$ which is homologous (in $Kh'$) to $\s_0$ and consists 
of terms with quantum grading $q>s-1$, which contradicts the equality $s([\s_0])=s-1$.  

In the case of $\ZZ/2$ coefficients, Lee's theory as above does not produce a spectral sequence (indeed, the resulting 
homology is isomorphic to $Kh_{\ZZ/2}(K)$). However,  a modification of Lee's construction \cite{Tur} works in this 
case: one considers a filtered theory with multiplication and comultiplication maps 
\begin{eqnarray*}
 m(\up\otimes \up) &=& \up  \qquad \qquad \qquad \qquad \Delta(\up) = \up\otimes \um + \um\otimes\up +\up\otimes\up \\ 
 m(\um\otimes \um)&=& \um  \qquad  \qquad  \qquad \qquad \Delta(\um) = \um\otimes \um  \\
 m(\up\otimes \um) &=& m(\up\otimes \um)= \um.   
\end{eqnarray*}
As explained in \cite{Tur}, Lee's arguments go through  to yield a spectral sequence whose $E^2$-term is 
the Khovanov homology with $\ZZ/2$-coefficients, and the $E^\infty$-term is  $\ZZ/2\oplus \ZZ/2$ when 
$K$ is a knot. When the knot $K$ is given by a braid, the canonical generators for this theory 
are given by two elements $(\um+\up) \otimes \um \otimes (\um+\up) \dots$
and $\um \otimes (\um+\up) \otimes \um \dots$ obtained by labeling the alternate components of the canonical 
resolution of the braid by  $(\um+\up)$ and $\um$. The transverse invariant $\psi$ is again the lowest
quantum degree part of the canonical generators. Moreover, a variant of the $s$-invariant can be defined in the same way, and by \cite{MTV} it takes the same values as the original Rasmussen's $s$, so our argument from the preceding paragraph still applies.  

It remains to deal with the reduced case. To obtain the reduced Khovanov complex, one places \ a marked point 
on the knot, forms the subcomplex $CKh_{-}(K)$ by labeling the the marked circle by $\um$ in every resolution, 
and considers the quotient complex $CKh(K)/CKh_{-}(K)$. For $\ZZ/2$ coefficients,  
the spectral  
sequence of \cite{Tur} works just as well in the reduced case. There is only one canonical generator $\s_o$ that survives;   
for a braid, it is given by the cycle
$(\um+\up)\otimes \um \otimes(\um+\up)\otimes\dots= \up\otimes \um \otimes(\um+\up)\otimes\dots$
with the label of $(\um+\up)$ on the marked circle of the oriented resolution of the braid. 
The lowest quantum degree part is $\up\otimes \um\otimes \dots \um$,  
which is precisely the reduced version of the transverse invariant.
The quantum grading shifts by 1 in the reduced case: we have 
$s^{red}(\s_o) = s_{min}^{red}=s(K)$, and the transverse invariant lives 
in the component $Kh^{0, sl+1}$ of the reduced homology.  This does not 
affect the validity of our argument.

To show the converse, note that for a $Kh_{\ZZ/2}$-thin link, 
reduced $Kh^{0, *}$ can only be non-trivial for one value 
of the quantum grading, namely $q=s$. For the reduced version of the transverse invariant, 
$q(\psi)=sl+1$, and the result follows. 

 (It is perhaps worth pointing out that the non-reduced 
Khovanov homology over $\ZZ/2$ is the direct sum of two copies of the reduced homology,
 $Kh_{\ZZ/2}^{n, q} = Kh_{\ZZ/2, red}^{n, q-1} \oplus Kh_{\ZZ/2, red}^{n, q+1}$ \cite{ORS}. Thus,
 if we define thin knots as those whose   
 non-reduced homology is supported on  two diagonals, or those   
 whose  reduced homology is supported on one diagonal, the set of thin 
 knots will be the same in both reduced and non-reduced cases.)  
\end{proof}

\begin{proof}[Proof of Corollary \ref{tightness}] We only need to recall two facts: by \cite{MO}, 
quasi-alternating knots are $Kh_{\ZZ/2}$-thin and have $s=\sigma$, and by \cite{ContOS}, non-vanishing 
of the contact invariant $c(\xi)$ is sufficient for tightness 
of a contact structure $\xi$.
\end{proof}

\section{Examples} \label{examples} 

In this section we give examples of contact structures whose tightness can be established by using Corollary 
\ref{tightness}. Since the result would be trivial for quasipositive braids, we are looking for non-quasipositive knots such that the $s$-bound for their self-linking number is nevertheless sharp. 
Among knots with 10 crossings or less, there are exactly three such knots, namely the mirrors 
of $10_{125}$, $10_{130}$ and $10_{141}$ in the Rolfsen table. (As was indicated to the author by Lenny Ng, this can be seen by contrasting the list of quasipositive knots from \cite{Baa} and the values of the maximal self-linking numbers \cite{Ng}.) We use each of these knots to obtain an infinite family of tight contact 
structures, and show that  contact structures in two of these families  are not Stein fillable.

\begin{example} \label{e125} 

For $r\geq 5$,  consider the pretzel link $P(-r, 3, -2)$ (for $n=5$, this is the mirror of the knot $10_{125}$ in the Rolfsen table), and let $K_r$
be its transverse representative given by the closed braid 
$$
(\sigma_1)^{-r} \sigma_2 \sigma_1^{3} \sigma_2.
$$ 
We can use the algorithm from \cite{HKP} to obtain the contact surgery description 
for the induced contact structure $\xi$ on the branched double cover  $\Sigma(K_r)$.
For $r=5$, we get the surgery diagram
shown on the left of Figure~\ref{125}; when $r>5$, we have the diagram  with $r$ $(+1)$-surgeries instead of five. 
(Strictly speaking, \cite{HKP} gives a slightly different surgery diagram shown on the right of Figure~\ref{125}.
The two unoriented  surgery links can be easily shown to be Legendrian isotopic, and we prefer the more symmetric 
diagram. In other examples below, we will also pick  surgery links slightly different from but Legendrian isotopic to 
those given by \cite{HKP}.)     

\begin{figure}[htb] 
\includegraphics[scale=1.0]{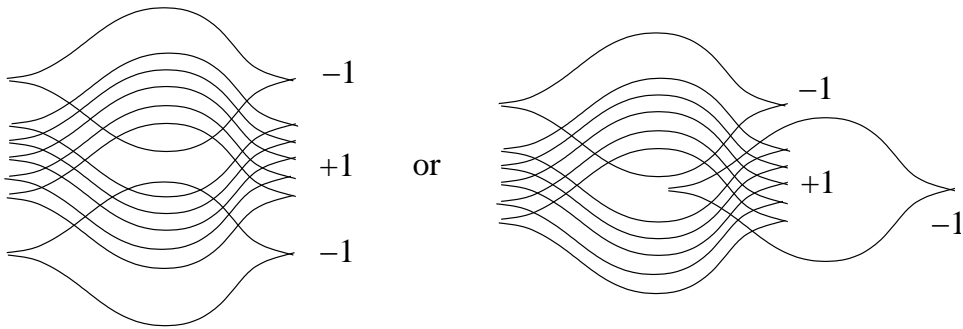}
\caption{The surgery diagram for the branched double cover of the transverse knot $K=K_5$.}\label{125} 
\end{figure}

The underlying smooth manifold  $\Sigma(K_r)$ is  the Seifert fibered space $M(-1; 2/3, 1/2, 1/n)$.  (See Figure~\ref{kirby125} for a sequence of Kirby calculus  moves demonstrating this for $n=5$.) 

Each link $K_r$ is quasi-alternating. Indeed, $|\det (K_r)| = |H_1 ( M(-1; 2/3, 1/2, 1/n) )| = r+6$. On the other hand, 
resolving the crossing circled in Figure \ref{braid125} in two possible ways, we obtain the link $K_{r-1}$ and the unknot. Repeating the procedure $r$ times, we get the link $K_0$, which is the trefoil linked once with the unknot.
Thus $K_0$ is an alternating link with $|\det (K_0)|=6$, and, since $|\det (K_r)| = |\det(K_{r-1})|  + |\det(\text{unknot})|$,  we see by induction that $K_r$ is quasi-alternating.

\begin{figure}[ht] 
\includegraphics[scale=0.5]{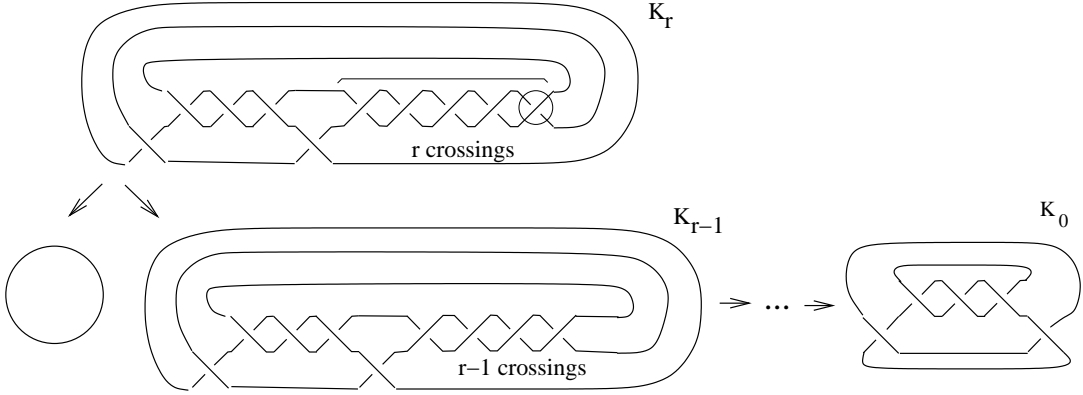}  
\caption{The links $K_r$ are quasi-alternating.}\label{braid125}
\end{figure} 

We next check the hypothesis of Theorem \ref{nonzero} when $r$ is odd (i.e. $K_r$ is a knot). We compute $sl(K)= 2-r$. 
The knot $K_r$ is $Kh$-thin, so $s$ equals to the signature $\sigma(K_r)=3-r$ (we compute the signature via 
the Goeritz matrix of the knot \cite{GL}).   
  
When $r$ is even, $K_r$ is a two-component link, so  Theorem \ref{nonzero} does not apply. However, we can argue that 
$\psi(K_r) \neq 0$ by \cite[Theorem 4]{Pla1}, since $\psi(K_{r-1}) \neq 0$, and the transverse braid $K_{r-1}$
is obtained from  $K_{r-1}$ by resolving a negative crossing. 

Theorem \ref{roberts} now implies that the branched double cover of each $K_r$ is a tight contact manifold. 
We now show that none of them are Stein fillable. Since $\Sigma(K_5)$ can be obtained from any of $\Sigma(K_r)$ by 
a sequence of Legendrian surgeries, it suffices to consider the contact structure $\xi_K = \xi$ that corresponds to 
$K=K_5$   and is shown on Figure~\ref{125}.

\begin{figure}[htb!] 
\includegraphics[scale=0.75]{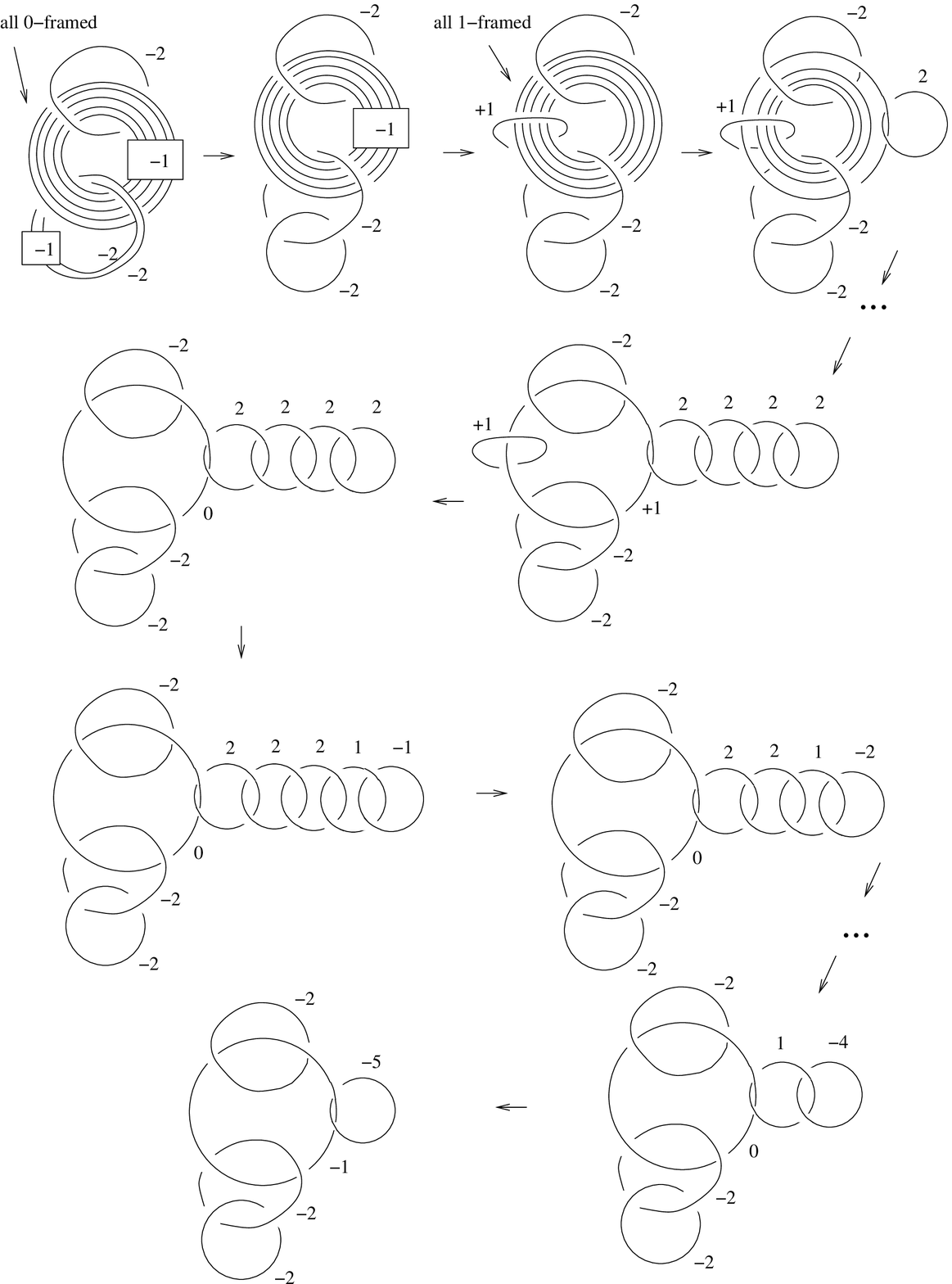}
\caption{Kirby moves.} \label{kirby125}
\end{figure}  

\begin{figure}[htb!] 
\includegraphics[scale=0.8]{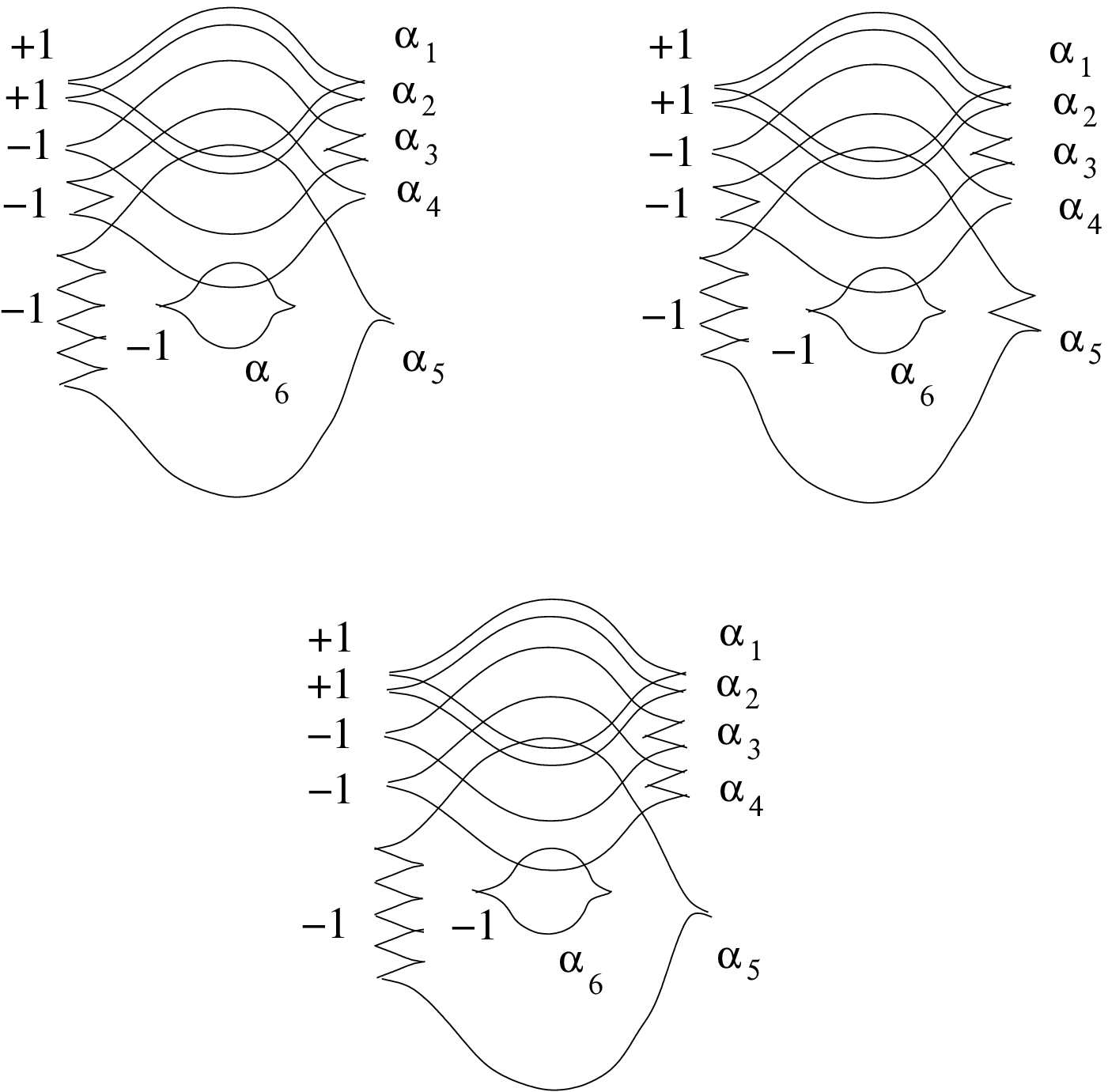}
\caption{The manifold $Y=M(-1; 2/3, 1/2, 1/5)$ carries three tight contact structures: $\xi_1$ (top left), $\xi_2$ (top right), and $\Xi$ (bottom).}\label{3cont} 
\end{figure}  

We have already mentioned that the branched double cover of $K$ is the Sefert fibered space $Y=M(-1; 2/3, 1/2, 1/5)$. 
Tight contact structures on this space were 
classified in \cite{GLS}; $Y$ carries three tight contact structures $\xi_1$, $\xi_2$ and $\Xi$ given by surgery diagrams on Figure~\ref{3cont}. To identify our contact structure $\xi$ among these three, we compute their $d_3$ invariants.

Recall \cite{DGS} that the three-dimensional invariant $d_3$ of a contact structure given by a contact  surgery
diagram can be computed as  
$$
d_3 (\xi) = \frac{c_1(\s)^2-2\chi(X)-3 \sign(X)+2}4 +m,  
$$  
where $X$ is a 4-manifold bounded by $Y$ and obtained by adding $2$-handles to $B^4$ as dictated by the surgery 
diagram, $\s$ is the corresponding $\Spinc$ structure on $X$, and $m$ is the number of $(+1)$-surgeries in the diagram.  
The $\Spinc$ structure $\s$ arises from an almost-complex structure defined in the complement of a finite set in $X$,
and the class  $c_1(\s)$ evaluates on each homology generator of $X$ corresponding to the handle attachment 
along an (oriented) Legendrian knot as the rotation number of the knot. 

For the contact structure  $\xi$  on $Y=M(-1; 2/3, 1/2, 1/5)$  defined by the surgery diagram from Figure~\ref{125}, 
we compute  $c_1(\s)=0$ and $d_3(\xi)=-\frac12$.

Let $\alpha$ be a Seifert surface a component of the Legendrian surgery link capped off 
by  the core of a handle attached along this component;
$c_1(\s)$ evaluates on $\alpha$ as the rotation number of the corresponding Legendrian knot. 
The classes of such surfaces $\alpha$ generate
$H_2(X)$; labelling the components of Legendrian surgery links on Figure~\ref{3cont} as shown, we compute the Poincar\'e duals:  
\begin{align*}
PD c_1(\xi_1)&= \frac1{11}\left(-29 \alpha_1  -29 \alpha_2 +20 \alpha_3 +12 \alpha_4  -3\alpha_5 +6 \alpha_6 
\right),  \\
PD c_1(\xi_2)&=\frac1{11}\left(-17 \alpha_1 - 17 \alpha_2 +14 \alpha_3 + 4 \alpha_4  - \alpha_5 -2 \alpha_6 
\right),   \\
PD c_1 (\Xi)&=  \alpha_1+\alpha_2-\alpha_5,  
\end{align*}

and thus 
$$
d_3(\xi_1) =\frac1{22} \quad, d_3(\xi_2) =\frac5{22} \quad,   d_3(\Xi) = -\frac12. 
$$
Because for the contact structure $\xi$ from Figure~\ref{125} we have 
$d_3(\xi)= -\frac12$,   it follows that 
$\xi$  is in fact the contact structure $\Xi$.

We show that $\xi$ is not Stein fillable, combining the ideas from \cite{GLS} and \cite{Li}. More precisely, 
we will show that $Y$ carries no Stein fillable contact structures with $d_3=-\frac12$.  
 
We first observe that $Y$ is an $L$-space, for example because it is a branched double cover 
of a quasi-alternating knot \cite{OS2}. It follows \cite{OS-L} that $b_2^+(X)=0$ for any symplectic filling $X$ of a contact 
structure on $Y$. By the argument in \cite{GLS}, this implies that $b_1(X)=0$. Now, observe that the space $-Y$
can be represented as the boundary of the plumbing shown on Figure \ref{-Yplumb}.
\begin{figure}[ht] 
\includegraphics[scale=0.75]{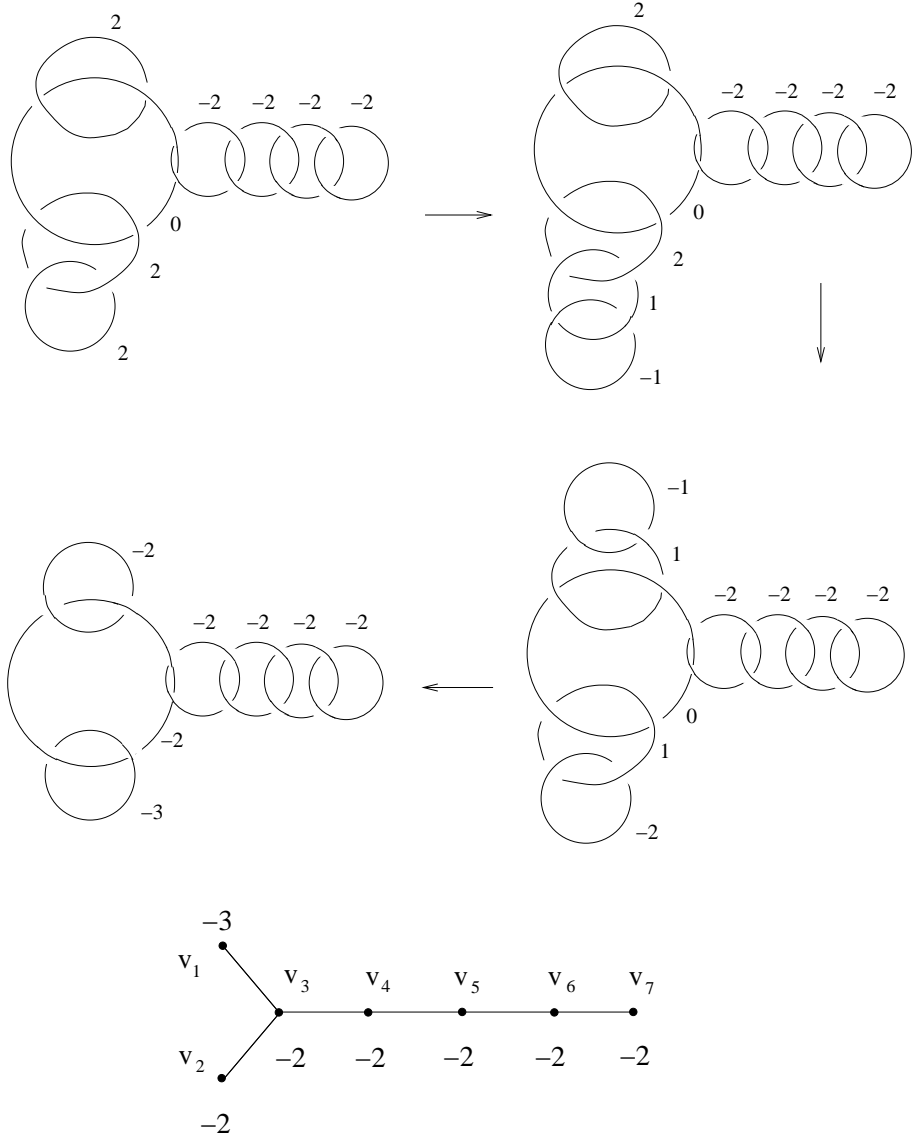}
\caption{Kirby diagrams for $-Y$ and the plumbing graph for $W$.} \label{-Yplumb}
\end{figure}

Denote by $W$ the 4-manifold with boundary $-Y$ given by this plumbing. If $X$ is a symplectic filling for $\xi$, 
then $X\cup W$ is an oriented negative-definite closed 4-manifold. By Donaldson's theorem, 
the intersection form on $X \cup W$ is standard diagonal $\langle -1 \rangle ^n$. To get restrictions on the intersection form 
of $X$, we consider the embeddings of the lattice given by Figure \ref{-Yplumb} into the standard negative-definite 
lattice, following \cite{Li}. Let $e_i$, $i=1, 2, \dots n$, be the basis of $ \langle -1 \rangle ^n$ such that $e_i \cdot e_j = -\delta_{ij}$.  
Let $v_i$ be the basis of  $H_2(W)$ corresponding to the vertices of the plumbing graph of Figure \ref{-Yplumb}. Up to permutations and sign reversals of $e_i$ (which are automorphisms of the lattice $ \langle -1 \rangle ^n$), we have 
\begin{align*}
&v_3 \mapsto e_1+e_2, \quad v_2 \mapsto -e_1 +e_3,  \quad v_1 \mapsto -e_1 -e_3 +e_5, \\     &v_4 \mapsto -e_2 +e_4, \quad
v_5 \mapsto -e_4 +e_6, \quad v_6 \mapsto -e_6+e_7, \quad v_7 \mapsto -e_7+e_8                                           
\end{align*}
(Another possibility would be for the first four vectors to embed as 
$$
v_3 \mapsto  e_1+e_2, \quad v_2 \mapsto -e_1 +e_3, \quad v_1 \mapsto -e_2 + e_4 + e_5, \quad v_4 \mapsto -e_1 - e_3, 
$$
but this leads to a contradiction when we try to embed $v_5$.)

The orthogonal complement $L$ of the image of the lattice generated by images of $v_i$'s in $\langle -1 \rangle ^n$ is then spanned by 
the vectors
$$
-e_1+e_2-e_3+e_4-2e_5+e_6+e_7+e_8, \quad
 e_9, \quad \dots \quad e_n,
$$
and the intersection form on $L$ is the diagonal form $ \langle -11 \rangle  \oplus \langle -1 \rangle ^{n-8}$.    
Because $H_1(Y) = \ZZ/11$ (indeed, $|H_1(Y)| = \det(10_{125})=11$), and both $H_2(X)$, $H_2(W)$ are torsion-free, we have 
$$
0 \to H_2(X)\oplus H_2(W) \to H_2(X\cup W) \to \ZZ/11 \to 0,  
$$  
and thus $H_2(X,\ZZ)$ is a subgroup of  $L= \ZZ^{n-7}$ of index $11$. Set $m=n-7=b_2(X)$, and let 
$\{u_1,\, u_2, \, \dots u_m\}$ be basis 
of $L$ in which the form is diagonal, and $u_1 \cdot u_1 = -11$. The vectors $11 u_1$,  $11 u_2$, ... $11 u_m$
lie in $H_2(X, \ZZ)$, and generate $H_2(X, \QQ)$ over $\QQ$.

Now, assume that $(X, J)$ is a Stein filling for $\xi$, 
and $\s_J$ is the corresponding 
$\Spinc$ structure on $X$. Let $\bar \xi$ be the contact structure on $Y$ conjugate
to $\xi$; then $\bar \xi$ has a Stein filling $(X, -J)$, with $\s_{-J}= \bar{\s}_J$ the corresponding 
$\Spinc$ structure. We have $d_3(\bar \xi)= -\frac{1}2$, and the classification of contact structures 
on $Y$ implies that $\bar \xi$ is isotopic to $\xi$. Then by \cite{LM} we must have $c_1(\s_J)= c_1(\s_{-J})$,
so $c_1(\s_J)=0$.

On the other hand, 
 $c_1(\s)$ evaluates as an odd integer on each vector   $11 u_1$, $11 u_2$, $\dots 11 u_m$; it follows that 
 $m=0$. Then $d_3(\xi)=0$, which contradicts the calculation $d_3(\xi)=-\frac12$.

\end{example} 

\begin{example}  \label{e141} Consider the transverse representative of the mirror of the knot $10_{141}$ given by the 
braid $\sigma_1^{-4} \sigma_2 \sigma_1^{3} \sigma_2^2$. We consider the family of braids  
$$
K_r = \sigma_1^{-r} \sigma_2 \sigma_1^{3} \sigma_2^2.
$$
The contact surgery description for the corresponding contact structures are shown on Figure \ref{141};
the surgery diagrams are quite similar to those in the previous example, but have one extra component.
The Kirby calculus moves similar to those in Figure \ref{kirby125} show that the branched double cover 
is the Seifert fibered space $M(-1; 2/3, 2/3, 1/n)$. 
\begin{figure}[htb] 
\includegraphics[scale=0.9]{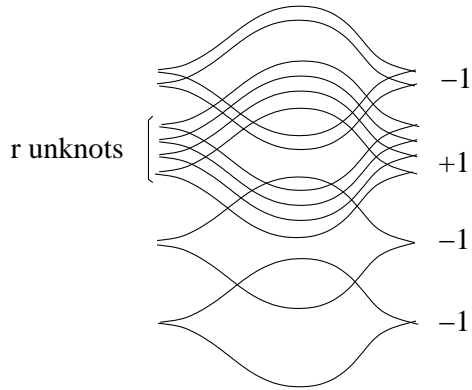}
\caption{The surgery diagrams for the branched double covers of the transverse links $K_r$.}\label{141} 
\end{figure}

As before, we can show that all the braids $K_r$ are quasi-alternating. Indeed, 
we resolve one of the negative crossings to obtain $K_{r-1}$ and a trefoil as two resolutions;
we also observe that $K_0$ is the connected sum of two trefoils.
Since  $|\det (\text{trefoil})| =3$, $|\det( K_0)|=9$ and$ |\det (K_r)| = |H_1 ( M(-1; 2/3, 1/2, 1/n) )| = 9+3r$,  
each $K_r$ is quasi-alternating by induction. 

Next, we compute $sl(K_r)=3-r$, and $s=\sigma(K_r)= 4-r$; the hypotheses of Corollary \ref{tightness} are therefore 
satisfied, and all branched covers $\Sigma(K_r)$ are tight contact manifolds.

For the contact structure on the branched cover of $K_4$, we compute $d_3=0$, which provides no obstruction to Stein fillability.  However, for the braid $K_6$ we get $d_3=-\frac12$. We then argue as in the previous example 
to show that the branched cover of $K_6$ is not Stein fillable 
(and thus the branched double covers of all braids  $K_r$ with $r\geq 6$ are not Stein fillable either). 
\begin{figure}[htb] 
\includegraphics[scale=0.85]{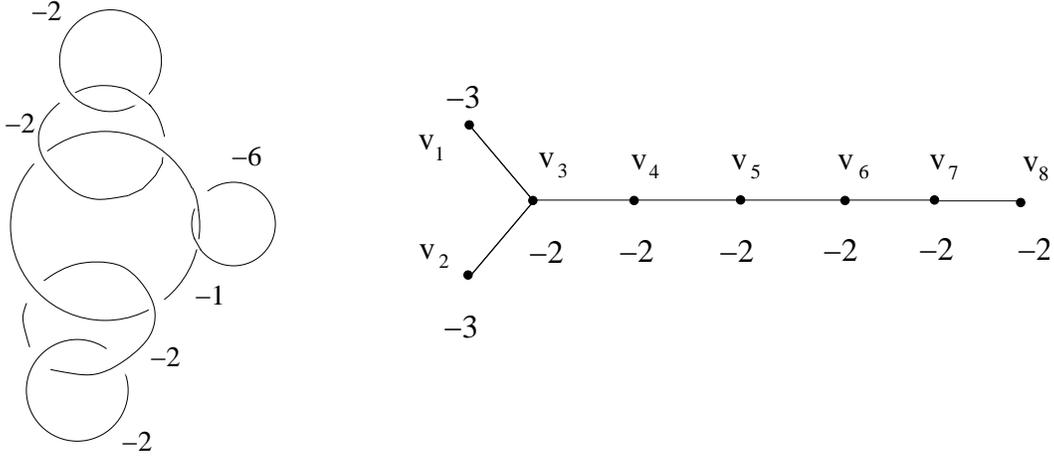}
\caption{A Kirby diagram for $Y= \Sigma(K_6)$ and the plumbing graph for $-Y$ of Example \ref{e141}.}\label{mY141} 
\end{figure}  
Denote $Y= \Sigma(K_6) = M(-1; 2/3, 1/2, 1/6)$; then $-Y$ is the boundary of the plumbing $W$ given by the graph 
on Figure \ref{mY141}.  As before, for any symplectic filling $X$ of $Y$  
the union  $X\cup W$ is a negative-definite closed 4-manifold with the  
standard diagonal intersection form. Up to changing the signs and the order of the vectors $e_i$ in the diagonal basis, there is a unique embedding
 of the lattice given by Figure \ref{mY141} into 
$\langle -1 \rangle ^n$, given by 
\begin{align*}
&v_3 \mapsto e_1+e_2, & &v_1 \mapsto -e_1 -e_3 +e_5, & &v_2 \mapsto -e_1 +e_3+e_4, & &v_4 \mapsto -e_2 +e_6,\\      
&v_5 \mapsto -e_6 +e_7, & &v_6 \mapsto -e_7+e_8,    & &v_7 \mapsto -e_8+e_9, & &v_8 \mapsto -e_9+e_{10},   
\end{align*}
and thus the orthogonal complement of this lattice in $\langle -1 \rangle ^n$ is 
$ \langle -9 \rangle  \langle -3 \rangle \oplus \langle -1 \rangle ^{n-10}$.
As in the previous example, the classification of tight contact structures on  $M(-1; 2/3, 1/2, 1/6)$
\cite{GLS} impies that our contact structure is isotopic to its conjugate, and so $c_1(X)=0$ for any Stein filling.  
Since $|H_1(Y)| = 27$, similar parity argument  
shows that $b_2(X)=0$, and so $d_3$ must be zero, a contradiction. 

\end{example}

\begin{remark} One can try to argue as in \cite{GLS} to investigate symplectic fillability in
Examples \ref{e125} and \ref{e141}: a slightly more involved
agrument modulo 8 puts further restrictions on the value $d_3$ for symplectic fillings 
(with diagonal odd intersection form). However, this gives no obstruction to symplectic fillability 
of any contact structures in the above two examples.

 

In the opposite direction, certain tight open books with the punctured torus page and pseudo-Anosov 
monodromy can be shown 
to be symplectically fillable as perturbations of taut foliations \cite{HKM1}.  
We note that our examples are not pseudo-Anosov, so these results do not apply.

\end{remark}

\begin{example} \label{e130}
A transverse representative of the mirror of $10_{130}$ with the maximal self-linking number 
is given by the braid $\sigma_1^{-3} \sigma_2 \sigma_1^2 \sigma_2^2 \sigma_3 \sigma_2^{-1} \sigma_3$.
We consider a family of transverse braids 
 $$
  K_r=\sigma_1^{-r} \sigma_2 \sigma_1^2 \sigma_2^2 \sigma_3 \sigma_2^{-1} \sigma_3.
  $$
First, we check that all the underlying links are quasi-alternating. Resolve of the negative crossings  
among those given by $\sigma_1^{-r}$ to obtain $K_{r-1}$ as one of the resolutions and the unknot as the other.
Observe that $K_0$ is a two-component alternating link of $\det =14$ (with $5_2$ knot and the unknot as 
components, linked once). Finally, compute $|\det (K_r)| = 14+r$ (one way to see this is to compute the size of $H_1$
of the branched double cover of $K_r$ which is a Seifert fibered space shown on Figure \ref{130}). 

\begin{figure}[htb] 
\includegraphics[scale=0.9]{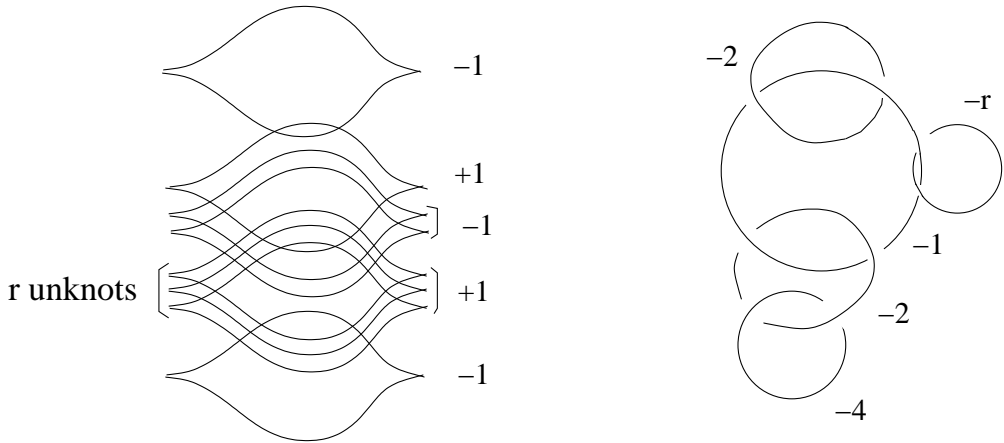}
\caption{The surgery diagrams for the branched double covers of the transverse links $K_r$.}\label{130} 
\end{figure}
  
The hypothesis of Corollary \ref{tightness} holds: $sl(K_r)= 2-r = \sigma-1$. Therefore, the branched double
covers of the transverse links $K_r$ are all tight.   
  
We do not investigate the fillability question in this case. (One can still try to use the classification 
of tight contact structures on these small Seifert fibered spaces, the fact they are all $L$-spaces, 
and non-vanishing of $d_3$ for some values of $r$, but the intersection form for the corresponding plumbings is harder to analyze.)  
\end{example} 

\begin{remark} In Examples \ref{e125} and \ref{e141}, transverse links are 3-braids, and 
the contact structures on the branched double covers can be given by open books whose page 
is a once-punctured torus. Tightness of these contact structures can be established by using results 
of \cite{HKM} or (easier yet) by rewriting the braids in the ``standard'' form and using Baldwin's work \cite{Ba}.
Example \ref{e130} deals with 4-braids; the page of the corresponding open books is a twice-punctured torus, 
and known results do not apply.  
\end{remark}

\begin{remark} In all of the above examples, we  checked explicitly that our families of links are quasi-alternating. 
In fact a weaker condition, $\rk Kh_{\ZZ/2}(K) = |\det (K)|$ is sufficient to ensure that the spectral 
sequence from $Kh$ to $\HF$ collapses at the $E^2$ stage. For any individual reasonably small knot this 
can be checked by a computer, for example using Baldwin's {\tt Kh} program \cite{Ba1} that computes the rank of reduced Khovanov 
homology with $\ZZ/2$ coefficients. Checking the second condition, $sl(K)=s-1$, is also routine for $Kh_{\ZZ/2}$-thin
knots (alternatively, one can use the {\tt Trans} program \cite{Ba1} to check $\psi\neq 0$). 
Thus tightness of the contact 
structure on the branched double cover can be established by a computer calculation.    
\end{remark}

\begin{example} \label{nex}
\begin{figure}[htb!] 
\includegraphics[scale=0.9]{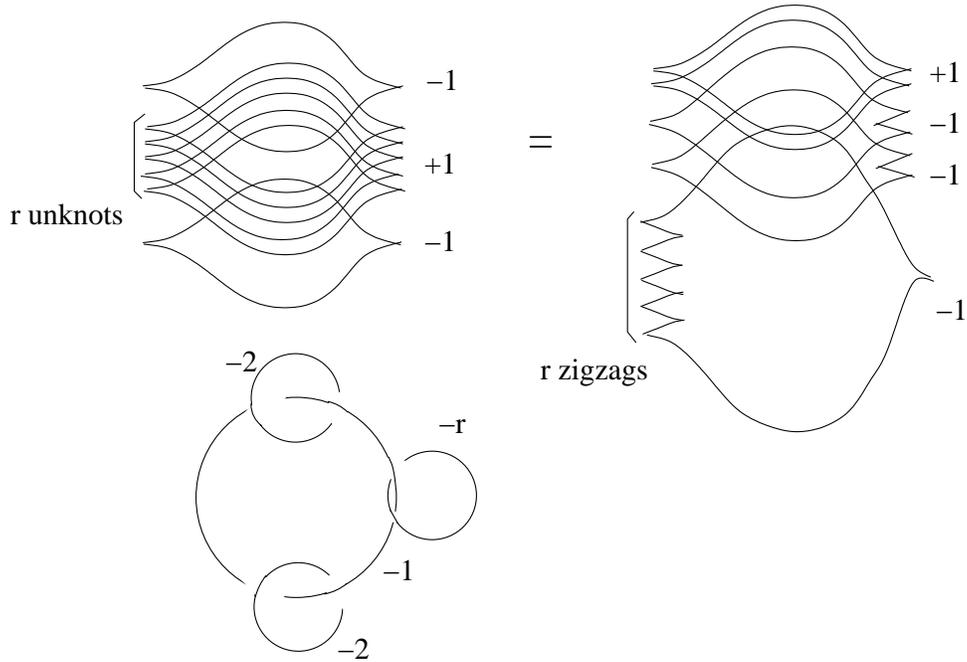}
\caption{The surgery diagrams for the branched double covers of the transverse links of Example \ref{nex}.}\label{non-ex} 
\end{figure}  
 John Baldwin has pointed out that for a family of (non-quasi-alternating) 
transverse 3-braids 
$$
K_r=\sigma_1^{-r} \sigma_2 \sigma_1^2  \sigma_2
$$
for $r>2$ the transverse invariant $\psi$ vanishes, but the contact structures on the branched 
double covers are tight, and have $c(\xi)\neq 0$. The vanishing of $\psi$ can be established by 
the computer program \cite{Ba1}, while non-vanishing of the contact invariant follows from calculations
in \cite{HKM} and \cite{Ba}. Obviously, Khovanov homology fails to detect tightness in this case 
(so this is really a non-example), but it is interesting to take a look at the corresponding contact 
structures. They are given by surgery diagrams on Figure \ref{non-ex}, and are very similar to 
the contact structures from Example \ref{e125}. As the latter are obtained from the former by Legendrian 
surgery on a knot, the  contact structures $\xi(K_r)$  cannot be Stein fillable. As 
the underlying smooth manifold is $M(-1;1/2, 1/2, 1/r)$ carries a unique Stein non-fillable contact structure 
for each $r$, these contact structures are precisely  those considered in \cite{GLS}, where 
most of them are shown to be symplectically non-fillable. One may wonder whether there is any relation between 
vanishing of $\psi$ and symplectic non-fillability (although such relation seems quite improbable). 
    
\end{example}


\begin{thebibliography}{HKM2}

\bibitem[APS]{APS} M.~Asaeda, J.~Przytycki, A.~Sikora, {\em Categorification of the Kauffman bracket skein module of $I$-bundles over surfaces}, Algebr. Geom. Topol. {\bf 4} (2004) 1177-1210.

\bibitem[Baa]{Baa} S. Baader, {\em Slice and Gordian numbers of track knots}, Osaka J. Math. {\bf 42} (2005), 257--271.

\bibitem[Bal1]{Ba} J. Baldwin, {\em Tight contact structures and genus one fibered knots},  Algebr. Geom. Topol.  {\bf 7}  (2007), 701--735. 

\bibitem[Bal2]{Ba1} J. Baldwin, {\tt Trans}, {\tt Kh} computer programs, available at 
{\tt http://math.columbia.edu/\~{}baldwin/} 

\bibitem[DGS]{DGS} F.~Ding, H.~Geiges and A.~Stipsicz, {\em Surgery diagrams for contact 3-manifolds},  Turkish J. Math.  {\bf 28},  no. 1, (2004),  41--74.

\bibitem[GLS]{GLS} P.~Ghiggini, P.~Lisca, A.~Stipsicz, {\em Tight contact structures on some small Seifert fibered 3-manifolds},  Amer. J. Math.  {\bf 129}  (2007),  no. 5, 1403--1447.

\bibitem[GL]{GL} C. Gordon, R. Litherland, {\em On the signature of a link},  Invent. Math.  {\bf 47}  (1978), no. 1, 53--69. 

\bibitem[HKM1]{HKM} K. Honda, W. Kazez, G. Mati\'c, {On the contact class in Heegaard Floer homology},
arxiv: math/0609734. 

\bibitem[HKM2]{HKM1} K. Honda, W. Kazez, G. Mati\'c,  {\em Right-veering diffeomorphisms of compact surfaces with boundary. II}, arxiv:math.GT/0603626.

\bibitem[HKP]{HKP}  S.~Harvey, K.~Kawamuro, O.~Plamenevskaya, {\em On transverse knots and branched covers},
arXiv: 0712.1557.

\bibitem[Kh]{Kh} M. Khovanov, {\em A categorification of the Jones polynomial},
 Duke Math. J.  {\bf 101}  (2000),  no. 3, 359--426.

\bibitem[Lee]{Lee} E. S. Lee, {\em An endomorphism of the Khovanov invariant},  
Adv. Math. {\bf 197} (2005), no. 2, 554--586. 

\bibitem[Li]{Li} P.~Lisca, {\em On symplectic fillings of $3$-manifolds}, Turkish J. Math. {\bf  23}  (1999),  no. 1, 151--159.

\bibitem[LM]{LM} P.~Lisca, G.~Mati\'c, {\em Tight contact structures 
and Seiberg-Witten invariants}, Invent. Math. {\bf 129} (1997), no. 3, 
509--525.

\bibitem[MO]{MO} C.~Manolescu, P.~Ozsv\'ath, {\em  On the Khovanov and knot Floer homologies of quasi-alternating links},
arXiv:0708.3249.  
 
\bibitem[MTV]{MTV} M. Mackaay, Paul Turner, Pedro Vaz, {\em A remark on Rasmussen's invariant of knots}, J. Knot Theory Ramifications  {\bf 16}  (2007),  no. 3, 333--344.

\bibitem[Ng]{Ng} L. Ng, {\em On arc index and maximal Thurston-Bennequin number},  arxiv:math/0612356. 
 
\bibitem[ORS]{ORS}  P.~Oszv\'ath, J.~Rasmussen, and Z.~Szab\'o, {\em Odd Khovanov homology}, arXiv:0710.4300. 
 
\bibitem[OS1]{OS-L}  P. Oszv\'ath and Z. Szab\'o, {\em Holomorphic disks and genus bounds},  Geom. Topol.  {\em 8},  (2004), 311--334. 
 
\bibitem[OS2]{ContOS} P. Oszv\'ath and Z. Szab\'o, {\em Heegaard Floer
  homologies and contact structures},  Duke Math. J. {\bf  129}  (2005),  no. 1, 39--61.

 
\bibitem [OS3]{OS2} P.~Ozsv\'ath,  Z.~Szab\'o, {\em On the Heegaard Floer homology
of branched double-covers}, Adv. Math.  {\bf 194}  (2005),  no. 1, 1--33.

 \bibitem[Pl1]{Pla1} O. Plamenevskaya, {\em Transverse knots and Khovanov homology}, {\em Math. Res. Lett.} \textbf{13}  (2006),  no. 4, 571--586. 

\bibitem[Pl2]{Pla2} O. Plamenevskaya, {\em Transverse knots, branched double covers and Heegaard Floer contact invariants}, {\em J.~Symplectic~Geom.}~\textbf{4}  (2006),  no. 2, 149--170.

\bibitem[Ra1]{Ra1} J. Rasmussen,
   {\em Floer homology and knot complements}, Ph.D. thesis, Harvard, 2003, ArXiv: math.GT/0306378.

\bibitem[Ra2]{Ra} J. Rasmussen,
   {\em Khovanov homology and the slice genus}, Invent. Math, to appear, ArXiv: math.GT/0402131.



\bibitem[Ro]{Ro} L. Roberts, {\em On knot Floer homology in double branched covers}, arXiv:0706.0741.
 
\bibitem[Sh]{Sh} A.~Shumakovitch, {\em Rasmussen invariant, Slice-Bennequin inequality, and sliceness of knots},
Arxiv: math.GT/0411643. 
 
\bibitem[Tu]{Tur} P. Turner, 
{\em Calculating Bar-Natan's characteristic two Khovanov homology},  J. Knot Theory Ramifications {\bf 15}
  (2006),  no. 10, 1335--1356.

 
\end{thebibliography}
\end{document}